\documentclass[a4paper,10pt]{article}
\usepackage[utf8x]{inputenc}

\usepackage{amssymb} 

\title{Monotonic lagrangian tori of standard and non standard types in toric and pseudotoric manifolds}
\author{Nikolai Tyurin\footnote{The author is partially supported by Laboratory of Mirror Symmetry NRU HSE, RF Government grant, ag. N 14.641.31.0001} \\  BLTPh JINR (Dubna) and NRU HSE (Moscow)}
\date{}

\begin{document}

\maketitle

\begin{abstract}

 In recent papers, summarized in survey [1], we construct a number of examples of non standard lagrangian tori  on compact toric varieties and as well on certain non toric
varieties which  admit pseudotoric structures. Using this pseudotoric technique we explain how non standard lagrangian tori of Chekanov type can be constructed and what is the topological 
difference  between standard Liouville tori and the non standard ones. However we have not discussed the natural question about the periods of the constructed twist tori; in particular 
the monotonicity problem for the monotonic case was not studied there. In the paper we present several remarks on these questions, in particular we show for the monotonic case
how to construct non standard lagrangian tori which satisify the monotonicity condition. First of all we study non standard tori which are Bohr - Sommerfeld with respect
to the anticanonical class. This notion was introduced in [2], where one
defines certain universal Maslov class for the ${\rm BS}_{can}$ lagrangian submanifolds in  compact simply connected monotonic symplectic manifolds.
Then we show how monotonic  non standard lagrangian tori of Chekanov type can be constructed. Furthemore we extend the consideration to pseudotoric setup
and construct examples of monotonic lagrangian tori in non toric monotonic manifolds: complex 4 - dimensional quadric and full flag variety $F^3$.

\end{abstract}

\section*{Introduction}

Our framework is as usual: let $(M, \omega)$ be a compact smooth simply connected real symplectic manifold of dimension $2n$  with integer symplectic form $\omega$, so $[\omega] \in H^2(M, \mathbb{Z})$.

Let $S \subset M$ be a compact orientable $n$ - dimensional submanifold. We say that $S$ is lagrangian if $\omega|_S  \equiv 0$, and Bohr - Sommerfeld of level $k \in \mathbb{Z}$ 
(or ${\rm BS}_k$ for short) if for any loop $\gamma \subset S$
and any disc $B_2 \subset M, \partial B_2 = \gamma$ one has $\int_{B_2} k \cdot \omega \in \mathbb{Z}$. Below we consider lagrangian tori only; however many things can be extended to much more general case.

  To establish for a given lagrangian torus $S \subset M$ is it Bohr - Sommerfeld or not one has to calculate the periods of $S$: take a basis $(a_1, ..., a_n) \in H_1(S, \mathbb{Z})$, realize it by
loops $\gamma_1, ...,  \gamma_n$, find discs, bounded by the loops, and compute their symplectic areas which gives the set of periods $(p_{a_1}, ...,  p_{a_n})$ defined up to $\mathbb{Z}$. Clearly
any lagrangian $S$ is ${\rm BS}_k$ if and only if the periods belong to $\frac{1}{k} \mathbb{Z}$; any other choice of the basis corresponds to a transformation of the period vector by
$SL(n, \mathbb{Z})$, therefore for any lagrangian torus one can define its Bohr - Sommerfeld level as the minimal $k$ such that every $k p_{a_i}$ belongs to $\mathbb{Z}$ if it exists
or if it does not saying that it is ${\rm BS}_{\infty}$. 

Fixing any almost complex structure $J$ on $M$, compatible with $\omega$, we get the complex determinant line bundle $K^{-1}_M = {\rm det} T^{1,0} M$ which we call anticanoncial line bundle
following the tradition. It depends on the choice of $J$, but its first Chern class does not being integer valued therefore it is a topological invariant of symplectic manifold.
Consider the case of monotonic symplecitc manifold namely when $c_1(K^{-1}_M) = k [\omega]$ for certain integer $k$. For this case we say that a lagrangian torus $S \subset M$ is
Bohr - Sommerfeld with respect to the anticanonical bundle (or ${\rm BS}_{can}$ for short) if it is Bohr - Sommmerfeld of level $k$.

{\bf Remark.} This assignment looks a bit artificial in general symplectic setup, but from the point of view of algebraic geometry it looks more natural. Indeed, consider
a Fano variety $X$. By the very definition its anticanonical line bundle $K^{-1}_X$ is ample therefore  certain power $(K^{-1}_X)^m$ induces an embedding of $X$ to the projective space.
Then choosing a standard Kahler form $\Omega$ on the projective space and restricting it to the image of $X$ we get a symplectic form $m \omega$ on $X$ such that single $\omega$
is integer as well since the class  $\Omega|_X = m c_1 (K^{-1}_X) = m [\omega]$. The resulting form $\omega$ can be called anti canonical; clearly it is not unique but
since all standard Kahler forms on the projective space are conjugated one could expect that the corresponding lagrangian geometries are the same. In this case ${\rm BS}_{can}$
lagrangian tori must characterize our given $X$ as well as algebro geometric ingredients do it. 

Suppose that $(X, \omega)$ is in addition monotonic, then for any ${\rm BS}_{can}$ lagrangian torus $S \subset M$ one has a universal class $m_S \in H^1(S, \mathbb{Z})$ defined in [2].
Fix a compatible almost complex structure $J$ on $M$ and an orientation on $S$. Then we get a realization of $K^{-1}_M$ together with the corresponding hermitian structure $h$
on it given by the hermitian triple $(G, J, \omega)$ (the Riemannian metric $G$ is reconstructed from two other elements) on the tangent bundle $T^{1,0}_J M$.  Then we choose
a hermitian connection $a \in {\cal A}_h(K^{-1}_M)$ whose curvature form $F_a$ is proportional to $\omega$. Then if $S \subset M$ is ${\rm BS}_{can}$ it implies that the restriction
$(K^{-1}_M, a)|_S$ admits a covariantly constant section $\sigma_S$ defined up to $\mathbb{C}^*$. On the other hand for a choosen orientation take the corresponding
section $\delta$ of the determinant ${\rm det} TS$, project it to the complex determinant $K^{-1}_M |_S$ and denote the result as $\delta_c$. Local computations ensure that
$\delta_c$ vanishes nowhere, therefore the ratio $\frac{\sigma_S}{\delta_c}: S \to \mathbb{C}^*$ defines an integer cohomology class given by the lifting the standard generator of $H^1(\mathbb{C}^*, \mathbb{Z})$,
and we denote this class as $m_S$. In [2] one shows that this class can be naturally understood as a universal version on the Maslov class; it is invariant under the Hamiltonian (isodrastic)
deformations of $S$. One can show that $S$ is monotonic if and only if $m_S$ is trivial, see [2]. 

Thus to study monotonic lagrangian tori we first search  ${\rm BS}_{can}$ lagrangian tori. Therefore first of all we study the periods.

\section{Toric case}

Now we come to our main subject: non standard lagrangian tori of Chekanov type. To start with let us take the simplest examples which have been  appeared many times. 

{\bf The first example.}  Consider the projective plane $M = \mathbb{C} \mathbb{P}^2$ with the standard symplectic form $\omega$. Fix homogenious coordinates $[z_0: z_1: z_2]$ and consider
 a pencil of plane conics $\alpha z_0 z_1 = \beta z_2^2$ with the base set $B = \{[0:1:0], [1:0:0] \}$.  The complement $\mathbb{C} \mathbb{P}^2 \backslash B$ is fibered over $\mathbb{C} \mathbb{P}^1$ with 
homogenious coordinates $[\alpha: \beta]$, and we denote this map as $\Psi$; all the fibers are smooth except two distinguished, when $\alpha$ or $\beta$ equals to 0. 
The real moment map 
$$
f_1 = \frac{\vert z_0\vert^2 - \vert z_1 \vert^2}{\sum \vert z_i \vert^2}
$$
 ``commutes'' with $\Psi$: its Hamiltonian action preserves the fibers of $\Psi$.
Therefore the data $(f_1, \Psi, B)$ define a pseudotoric structure on $\mathbb{C} \mathbb{P}^2$, see [1].

Then as it was shown in [1], any choice of a smooth loop $\gamma \subset (\mathbb{C} \mathbb{P}^1 \backslash \{[1:0], [0: 1]\})$ gives a lagrangian torus 
$$
T(\gamma, 0) = \bigcup_{p \in \gamma} (\Psi^{-1}(p) \cap \{f_1 = 0 \} ) \subset \mathbb{C} \mathbb{P}^2;
$$ 
if $\gamma$ is non contractible then $T(\gamma, 0)$ is of the standard type being Hamiltonian isotopic to a standard Liouville torus given by the toric structure on $\mathbb{C} \mathbb{P}^2$,
otherwise $T(\gamma, 0)$ is non standard lagrangian torus of Chekanov type. The proofs and discussion can be found in [1].

The projective plane is monotonic (being Fano), namely $K^{-1}_{\mathbb{C} \mathbb{P}^2} = {\cal O}(3)$ therefore the monotonicity coefficient equals $k = 3$ for this case.

Here we are interested in the Chekanov type tori, and the main claim is the following: it exists a smooth contractible loop $\gamma \subset  (\mathbb{C} \mathbb{P}^1 \backslash \{[1:0], [0: 1]\})$
such that the corresponding lagrangian torus $T(\gamma, 0)$ is monotonic. 

The construction is rather explicit: first we construct $T(\gamma, 0)$ which is ${\rm BS}_{can}$ and then we show that it is monotonic.
So we start with the periods. To calculate periods of $T(\gamma, 0)$ we need a ``section'' of the map $\Psi$.  Take projective line $\Sigma = \{ z_0 = z_1 \} \subset \mathbb{C} \mathbb{P}^2$,
and note that $\Sigma$ satisfies $f|_{\Sigma} \equiv 0$; on the other hand $\Sigma$ intersects each regular fiber of $\Psi$ exactly at two conjugated points. 
Then fix a smooth loop $\tilde \gamma_{Ch} \subset \Sigma$ such that: 1) the loop $\tilde \gamma_{Ch}$ does not intersect its
 image under the involution $\tau: [z_0: z_1: z_2] \mapsto [-z_0: - z_1: z_2]$
of our projective line $\Sigma$; 2) the loop $\tilde \gamma_{Ch}$ bounds a disc $B_2 \subset \Sigma$ of symplectic area $\frac{1}{3}$. Then for any point $ p \in \tilde \gamma_{Ch}$
with coordinates $[z_0(p): z_1(p): z_2(p)]$ take the loop $s_p = [z_0(p) e^{i t}: z_1 (p) e^{- it}: z_2 (p)]$. Note that $s_p$ never meets $\Sigma$ except at $p$ again (the condition 1) above),
therefore after the globalization $T = \bigcup_{p \in \tilde \gamma_{Ch}} s_p$ we get a smooth 2 -torus. The claim is that 1) this torus is lagrangian, non standard of Chekanov type;
2) it is ${\rm BS}_{can}$; 3) it is monotonic.

The first  follows from the fact that $\Psi$ restricted to $\Sigma$ is a double covering ramified exactly at $[1:0], [0:1] \in \mathbb{C} \mathbb{P}^1$ (and $\tau$
is exactly the conjugation of the fibers), thus $T = T(\gamma, 0)$ where $\gamma = \Psi(\tilde \gamma_{Ch}) \subset (\mathbb{C} \mathbb{P}^1 \backslash [1:0], [0:1])$.     
  On the other hand choose the classes of $\tilde \gamma_{Ch}$ and $s_p$ for any $p \in \tilde \gamma_{Ch}$ as the basis in $H^1(T, \mathbb{Z})$. Then the periods equal $(\frac{1}{3}, 1)$
since the choice of $\tilde \gamma_{Ch}$ and the fact that $s_p$ is an equatorial loop on conic $q_p = \{z^2_2(p) z_0 z_1 = z^2_0(p) z_2^2 \}$.

For the monotonicity we need some arguments from the toric geometry. Add the second moment map 
$$
f_2 = \frac{\vert z_2 \vert^2}{\sum_{i=0}^2 \vert z_i \vert^2}
$$
 (note that this function
preserves our section $\Sigma$ by the Hamiltonian action), then
the configuration of three projective lines $L_0 \cup L_1 \cup \subset \mathbb{C} \mathbb{P}^2, L_i = \{ z_i = 0 \}$, is invariant under the toric action generated by
moment maps $f_1$ and $f_2$. Therefore it exists a holomorphic section of the anticanonical bundle $\alpha_b \in H^0 (\mathbb{C} \mathbb{P}^2, K^{-1}_{\mathbb{C} \mathbb{P}^2})$
invariant under the full toric action. For our lagrangian torus $T(\gamma_{Ch}, 0)$ we fix a section of the determinant bundle ${\rm det} T(\gamma_{Ch}, 0)$
given by the following condition: take $\delta = X_{f_1}|_{T(\gamma_{Ch}, 0)} \wedge \phi^t_{X_{f_1}}(\nu)$ where $\nu$ is a non vanishing section of $T \tilde \gamma_{Ch} \subset T \Sigma$, 
lifted by the flow generated by $X_{f_1}$ on whole $T(\gamma_{Ch}, 0)$. Then the monotonicity condition for our non standard lagrangian torus $T(\gamma_{Ch}, 0)$
is equivalent to the following condition: non vanishing complex function $\alpha_b^*(\delta): \tilde \gamma_{Ch} \to \mathbb{C}^*$ has degree 1 (here $\alpha^*_b$ is top holomorphic form
with pole along $L_0 \cup L_1 \cup L_2$, dual to $\alpha_b$). 

Therefore our task is to calculate the degree of this map.

At each point of $\Sigma \backslash ([1:1:0], [0: 1:1])$ one has a tangent vector given by $X_{f_2}|_{\Sigma}$ since $f_2$ induces the slicing of $\Sigma$, corresponding
to the standard toric fibration of the projective plane. Without loss of generality we can suggest that $\tilde \gamma_{Ch}$
transversally intersects  the slices except at two tangent points.  Toric arguments imply that the degree $\alpha_b(\delta_0): \tilde \gamma_{Ch} \to \mathbb{C}^*$
is zero,  if  we take $\delta_0 = (X_{f_1} \wedge X_{f_2})|_{\tilde \gamma_{Ch}}$. But  we can express at each point $p \in \tilde \gamma_{Ch}$ tangent vector $\nu \in T_p \Sigma$
in terms of $X_{f_2}(p) \in T_p \Sigma$ using the complex structure $I$ (since $\Sigma$ is complex submanifold) namely it is a parametrization $\tilde \gamma_{Ch} (t), t\in [0; 2 \pi]$
such that $\nu (p(t)) = {\rm cos} (t) X_{f_2}(p(t)) + {\rm sin} (t) {\rm I}(X_{f_2}(p(t)))$; substituting this exrpession to $\alpha_b^*$ and using $\mathbb{C}$ - linearity
of the pairing we get that $\alpha_b^*(\delta): \tilde \gamma_{Ch} \to \mathbb{C}^*$ has degree 1. Therefore $T(\gamma_{Ch}, 0)$ is monotonic.

Note that this construction can not  directly lead to  construction of non standard  lagrangian torus of Chekanov type
of Bohr - Sommerfeld level 2: it does not exist a smooth equatorial loop $\tilde \gamma$ with trivial
intersection with $\tau(\tilde \gamma)$.

{\bf The second example.} Consider the direct product $ M = \mathbb{C} \mathbb{P}^1 \times \mathbb{C} \mathbb{P}^1$ endowed with the product symplectic form $\omega$ of type (1,1).
Thus $M$ is a complex 2 - dimnesional quadric.  Fix homogenious coordinates
on  both the projective lines $[x_0: x_1], [y_0: y_1]$ and consider pseudotoric structure $(f_1, \Psi)$ where 
$$
f_1 = \frac{\vert x_0 \vert^2}{\vert x_0 \vert^2 + \vert x_1 \vert^2} - \frac{\vert y_0 \vert^2}{\vert y_0 \vert^2 + \vert y_1 \vert^2}
$$
 and $\Psi: M \backslash B \to \mathbb{C} \mathbb{P}^1_w $ is given by $w_i = x_i y_i$. Here $[w_0: w_1]$ are the coordinates
on the last $\mathbb{C} \mathbb{P}^1_w$, and the base set $B$ is formed by 4 points which are the pairs of poles on the product projective lines.
In other words, $\Psi$ corresponds to pencil $\{ w_1 x_0 y_0 = w_0 x_1 y_1 \}$ of plane conics in $M$.  

A non standard lagrangian torus of Chekanov type arises from the following data: take a contractible loop $\gamma \subset \mathbb{C} \mathbb{P}^1_w \backslash \{[1:0], [0:1] \}$
then 
$$
T(\gamma, 0) = \bigcup_{p \in \gamma} \Psi^{-1}(p) \cap \{f_1 = 0 \} \subset M
$$
is a lagrangian torus. Again, as in the first example,  to compute the periods we need to find an appropriate ``section'' for the map $\Psi$ in $M$. In this case
take  rational curve $\Sigma = \{x_0 y_1 = x_1 y_0 \} \subset M$.  It admits involution $\tau: ([x_0: x_1] \times [y_0: y_1]) \mapsto ([-x_0: x_1] \times [-y_0: y_1])$.
The total symplectic area $\int_{\Sigma} \omega = 2$, therefore it is possible to find a smooth loop $\tilde \gamma_{Ch} \subset \Sigma$ such that 1) $\tilde \gamma_{Ch}$
does not intersect $\tau (\tilde \gamma_{Ch})$ in $\Sigma$; 2) the symplectic area of the small disc bounded by $\tilde \gamma_{Ch}$ equals $\frac{1}{2}$. Then we claim
that non standard lagrangian torus $T(\Psi(\tilde \gamma_{Ch}), 0) \subset M$ is ${\rm BS}_{can}$ and monotonic. Indeed, for this case the anticanonical class for $M$ has the type (2,2),
therefore ${\rm BS}_{can} = {\rm BS}_2$; at the same time periods for $T(\Psi(\tilde \gamma_{Ch}), 0)$ equal $(\frac{1}{2}, 1)$ since the loop
$l_p = \Psi^{-1}(p) \cap \{f_1 = 0 \}$ is equatorial for conic $\bar \Psi^{-1}(p) \subset M$ if $p$ is not a pole, and the total symplectic area of the conic equals to 2.

To establish monotonicity of $T(\gamma_{Ch}, 0) \subset M$ we use absolutely the same arguments as in the previous example. Adding the second moment map 
$$
f_2 = \frac{\vert x_0 \vert^2}{\vert x_0 \vert^2 + \vert x_1 \vert^2} + \frac{\vert y_0 \vert^2}{\vert y_0 \vert^2 + \vert y_1 \vert^2}
$$
 we consider a destinguished
section of the anticanonical line bundle $\alpha_b \in H^0(Q, K^{-1}_Q)$ which vanishes on the boundary divisor $D_b = \{ X_{f_1} \wedge X_{f_2} = 0 \} \subset Q$
formed by 4 lines; as in the previous case this section is invariant under the toric action generated by the moment maps $f_1, f_2$. The resting arguments are the same as above;
they lead to the main claim --- non standard lagrangian torus $T(\gamma_{Ch}, 0) \subset Q$ is monotonic. 

Essentially the same construction can be done for the case $M = \mathbb{C} \mathbb{P}^1 \times ... \times \mathbb{C} \mathbb{P}^1$ for $n$ copies of $\mathbb{C} \mathbb{P}^1$ equipped with the standard
product symplectic form $\omega = \oplus_{i=1}^n p_i^* \omega_{FS}$. Fix
homogenious coordinates $[x_i: y_i]$ for $i$ -th component, then the standard toric structure is given by the set of moment maps $(\tilde f_1, ..., \tilde f_n)$ 
where $\tilde f_i = \frac{\vert x_i \vert^2}{\vert x_i \vert^2 + \vert y_i \vert^2}$, the boundary divisor $D_b = \{X_{f_1} \wedge ... \wedge X_{f_n} = 0 \}$ consists of    
$2n$ components $D^1_i = \{ x_i = 0 \}, D_i^2 = \{ y_i = 0 \}$; the standard toric fibration admits monontonic lagrangian torus $T^n_{Cl} = \{\tilde f_1 = ... = \tilde f_n = \frac{1}{2} \}
\subset M$. 

In [1] one introduces a natural pseudotoric structure on $M$ given by the data $(f_1, ..., f_{n-1}, \Psi, B)$ where $f_i = \tilde f_i - \tilde f_{i+1}$, $\Psi: M \backslash B \to \mathbb{C} \mathbb{P}^1_w$
defined by $w_0 = x_1 ... x_n, w_1 = y_1 ... y_n$ and the base set $B$ consists of intersections $D_i^1 \cap D_j^2, i \neq j$. 

The values of the reduced moment maps $f_1, ..., f_{n-1}$ for the monotonic lagrangian torus $T^n_{Cl}$ equal 0, and therefore first of all we are looking for an approriate section
$\Sigma \subset M$ of the map $\Psi$ satisfies $f_i|_{\Sigma} = 0$.   It is not hard to see that the diagonal $\Sigma = \Delta = \{x_1 = ... = x_n; y_1 = ... = y_n \}$
is the desired section. The map $\Psi: \Sigma \to \mathbb{C} \mathbb{P}^1_w$ is the standard $n$ -covering; the symplectic area of $\Sigma$ equals $n$.

In  terms of the pseudotoric structure the standard lagrangian torus $T^n_{Cl}$ is given by the following choice: take smooth loop $\gamma_{Cl} = \{ [1: e^{it}] \} \subset \mathbb{C} \mathbb{P}^1_w$,
then $T(\gamma_{Cl}, 0)$ is exactly $T^n_{Cl}$. 

Now since $D_b$ represents the anticanonical class $K^{-1}_M = {\cal O}(2, ..., 2)$ while $[\omega] = c_1 ({\cal O}(1, ..., 1)) \in H^2(M, \mathbb{Z})$ the monotonicity coefficient
$k = 2$. Therefore to find a ${\rm BS}_{can}$ non standard lagrangian torus of Chekanov type we need a smooth constractible loop $\tilde \gamma_{Ch} \subset \Sigma \backslash (p_N, p_S)$
such that 1) after the rotation $\tau^i (\tilde \gamma_{Ch}) \cap \tilde \gamma_{Ch} = \emptyset$ where $\tau$ is the primitive rotation $e^{\frac{2 \pi i}{n}}$ of $\Sigma$;
2) the symplectic area of the disc bounded by $\tilde \gamma_{Ch}$ equals $\frac{1}{2}$ or equivalently $\frac{1}{2n}$ of the symplectic area of $\Sigma$. It is not hard to see
that such a loop exists. Then the corresponding non standard lagrangian torus $T(\gamma_{Ch}, 0)$ where $\gamma_{Ch} = \Psi(\tilde \gamma_{Ch}) \subset \mathbb{C} \mathbb{P}^1_w$,
is smooth and ${\rm BS}_{can}$.

The proof that $T(\gamma_{Ch}, 0)$ is monotonic follows the same scheme as in the first example. Essentially we have to check only one loop $\tilde \gamma_{Ch} \subset T(\gamma_{Ch}, 0)$.

Take the section $\alpha_b \in H^0(M, K^{-1}_M)$ with zeroset at $D_b$; this section is invariant under the Hamiltonian action of each $X_{f_i}$, where we add $n$-th moment map
$f_n = \sum_{i=1}^n \tilde f_n$; note that $X_{f_n}$ is tangent to $\Sigma$ at each point. Therefore non vanishing complex function
 $\alpha^*_b(X_{\tilde f_1} \wedge ... \wedge X_{\tilde f_n})|_{\tilde \gamma_{Ch}} : \tilde \gamma_{Ch}  \to \mathbb{C}^*$ has degree zero; here as above $\alpha^*_b$ is 
the top holomorhic form with pole at $D_b$ dual to $\alpha_b$. Again we can take a parametrization
$\theta: [0; 2\pi] \to \tilde \gamma_{Ch}$ such that the tangent vector $\frac{\partial}{\partial \theta} (p) = {\rm cos}(\theta) X_{f_n}(p) + {\rm sin} (\theta) {\rm I} (X_{f_n}(p))$
(without loos of generality we can take $\tilde \gamma_{Ch}$ convex). Consequently the degree of the map $\alpha_b(X_{f_1} \wedge ... \wedge X_{f_{n-1}} \wedge \frac{\partial}{\partial \theta})$
along $\tilde \gamma_{Ch}$ equals to 1; the section $\alpha_b|_{B_2}$ does not vanish being restricted to small disc $B_2 \subset \Sigma$ bounded by $\tilde \gamma_{Ch}$ therefore
the universal Maslov class must have trivial value on $[\tilde \gamma_{Ch}] \in H_1(T(\gamma_{Ch}, 0))$. Therefore $T(\gamma_{Ch}, 0)$ is monotonic.

  Consider now  general situation. Let $(M, \omega, \tilde f_1, ...., \tilde f_n)$ be a compact smooth simply connected toric symplectic manifold of dimension $2n$ with integer symplectic form
and complete set of commuting moment maps (first integrals) (see f.e. [3]). As it was shown in [1] it admits
a pseudotoric structure $(f_1, ..., f_{n-1}, \Psi, B)$ where $f_1, ..., f_{n-1}$ are pairwise commuting  moment maps derived from the complete set of first integrals by linear transformations, 
$B \subset M$ is the base set,
$\Psi: M \backslash B \to \mathbb{C} \mathbb{P}^1_w$ is the map with symplectic fibers, preserved by the Hamiltonian action of each $X_{f_i}$ (the ``commutation'' relation
for $\Psi$ and $f_i$). Recall the main idea of the construction.

For a given toric manifold $M$ we can fix the complex structure $I$ which is invariant with respect to the toric action, this structure is essentially unique. Consider the boundary 
divisor $D_b = \{X_{\tilde f_1} \wedge ... \wedge X_{\tilde f_n} = 0 \} \subset M$, --- it is a reducible complex subvariety of complex dimension $n-1$, formed by irreducible components
$D_1, ..., D_s$ where $s > n$. Taking the corresponding classes $[D_i] \in H_{2n-2}(M, \mathbb{Z})$ one can find a relation with integer coeficients  $\sum_{i=1}^s \lambda_i [D_i] = 0,
\lambda_i \in \mathbb{Z}$, coming from the combinatorial description of $M$. Recall that our toric $M$ can be recontructed from the corresponding convex polytop $P_M = F_{act}(M) \subset \mathbb{R}^n$ where
$F_{act} =  (\tilde f_1, ..., \tilde f_n)$ is the action map. Then for a given relation $\sum_{i=1}^s \lambda_i [D_i] = 0$ we separate the summands with non negative and negative coefficients
correspondingly, so $[D^+] = \sum_{\lambda_i \geq 0} \lambda_i [D_i]$, $[D^-] = \sum_{\lambda_j < 0} \lambda_j [D_j]$. Therefore we have $[D^+] = [D^-]$. Then it exists a holomorphic line bundle $L \to M$  which
admits  holomorphic sections $\alpha_{\pm} \in H^0(M, L)$ with zerosets $(\alpha_{\pm})_0 = D^{\pm}$. Since both $D^+$ and $D^-$ are invariant with respect to the toric action
generated by $X_{\tilde f_1}, ..., X_{\tilde f_n}$ it is a linear condition on $\tilde f_1, ..., \tilde f_n$ deriving a couple $f_1, ..., f_{n-1}$ of  moment maps which satisfy
the following condition: every element from the linear span $<D^+, D^-> \subset \vert L \vert$ is invariant under the action of each $X_{f_i}$. Thus a pseudotoric structure on $M$
is given by the data: the base set $B = D^+ \cap D^-$ is the base set of the pencil $<D^+, D^->$, the map $\Psi: M \backslash B \to \mathbb{C} \mathbb{P}^1_w$ is given by the elements of
the pencil, and the moment maps $f_1, ..., f_{n-1}$ are as described above.    

For the same $M$ probably there are many  relations of the form $\sum_i \lambda [D_i]=0 $, and formally we can proceed for any such relation; however as it was pointed out in [1]
  interesting and elegant examples appear in the ``middle'' cases when  $D^+$ is not too small or too big. Below we will see what would be a natural condition on $D^+$
corresponds to the ``middle'' case.  

Continue the construction. All the fibers $\Psi^{-1}(p)$ are smooth except for $p = [1:0], [0:1] \in \mathbb{C} \mathbb{P}^1_w$; the completions 
$$ D^+ = \overline{\Psi^{-1}([1:0])}, \quad  D^- = \overline{\Psi^{-1}([0:1])}  \subset M
$$
have as the supports the components of the  boundary divisor $D_b$. Any smooth function $h \in C^{\infty}(\mathbb{C} \mathbb{P}^1, \mathbb{R})$ can complete
the set of moment maps $f_1, ..., f_{n-1}$ --- take the lift $\Psi^* h$, correctly defined on the complement $M \backslash B$, then the last function must commute with every $f_i$. However
this lift can not be extended to whole $M$; at the same time every standard lagrangian torus from the toric fibration on $M \backslash D_b$ given by the moment maps $\tilde f_i$ is presented in the
 form $T(\gamma, c_1, ..., c_{n-1})\subset M$ where $\gamma \subset \mathbb{C} \mathbb{P}^1_w \backslash ([1:0], [0:1]) $ is a smooth non contractible loop  and $c_1, ..., c_{n-1}$ are
fixed values of the moment maps $f_1, ..., f_{n-1}$. 

Recall that in this setup the choice of a contractible smooth loop     $\gamma_{Ch} \subset (\mathbb{C} \mathbb{P}^1_w \backslash [1:0], [0:1])$ together with the choice of
the values give a smooth lagrangian torus $T(\gamma_{Ch}, c_1, ..., c_{n-1}) \subset M$ which is called non standard lagrangian torus of Chekanov type. The details can be found
in [1].

Suppose additionally that $(M, \omega)$ is monotonic, thus ${\rm P.D.} [D_b] = k [\omega] \in H^2(M, \mathbb{Z})$,  and suppose that  standard torus $T_{Cl} = \{\tilde f_1 = ...= \tilde  f_n = 0 \}$ is monotonic.
 Then we would like to find a  smooth ${\rm BS}_{can}$ non standard lagrangian torus of Chekanov type $T(\gamma_{Ch}, 0, ..., 0) \subset M$ (or just $T(\gamma, 0)$ for short) and establish whether or not it is 
monotonic.    

It is clear that the periods for lagrangian tori $T_{Cl}$ and $T(\gamma_{Ch}, 0)$ are almost the same except for one cycle in $H^1(T, \mathbb{Z})$ projected by $\Psi$
either to non contractible  loop $\gamma_{Cl}$ or to contractible $\gamma_{Ch}$. Thus we have to find a contractible loop with the desired  period.

  Following the same strategy as in the examples above, first we need an appropriate ``section'' of the map $\Psi: M \backslash B \to \mathbb{C} \mathbb{P}^1_w$, so a symplectic Riemann surface $\Sigma \subset M$
such that $\Sigma$ intersects each fiber $\Psi^{-1}(p)$ at finite number of points. Consider  the subset $N_0 = \{f_1 = ... = f_{n-1} = 0 \} \subset M$; there one has
a complex one - dimensional distribution $\pi$ given by $\pi(x)  = T_x N_0 \cap I(T_x N_0) \subset T_x M$ where $I$ is our toric complex structure; this distribution is invariant under 
the action of $T^n$ therefore it is integrable.
Indeed, for every $\tilde f_i$ from the  complete set $(\tilde f_1, ..., \tilde f_n)$ the Hamiltonian action $\phi^t_{X_{f_i}}$ preserves all the data which generate the distribution. 

The leaves of the integrable distribution are interchanged by the Hamitonian action of $f_1, ..., f_{n-1}$ and we can find the last $f_n$ as a linear combination of $\tilde f_1, ...  \tilde f_n$ such that
$X_{f_n}$ is tangent to the leaves; since for each leaf this action has exactly two fixed points, in general the leaves are rational curves. 

Fix a leaf of the distribution $\pi$ and denote it as $\Sigma \subset M$. By the  construction $\Sigma$ intersects each fiber $\Psi^{-1}(p)$ at a finite set of points (algebraic geometry provides the shortest way to see this fact), therefore the restriction $\Psi: \Sigma \to \mathbb{C} \mathbb{P}^1$ gives a finite covering ramified exactly at two points $[1:0]$ and $[0:1]$. Denote the degree of the covering
as $d$. Let the symplectic area of $\Sigma$ equals to $m$. What we need is a smooth contractible loop $\tilde \gamma_{Ch} \subset \Sigma \backslash (p_N, p_S)$, where $p_N, p_S$
are the branching points, such that two properties hold: 1)$\Psi(\tilde \gamma_{Ch}) \subset \mathbb{C} \mathbb{P}^1_w \backslash ([1:0], [0:1])$ is smooth , 2) $\tilde \gamma_{Ch}$ cuts a disc on $\Sigma$
of symplectic area $\frac{l}{k}$ where $l \in \mathbb{Z}$. It is clear that given such $\tilde \gamma_{Ch}$ the corresponding non standard lagrangian torus $T(\Psi(\tilde \gamma_{Ch}), 0)$
must be ${\rm BS}_{can}$.

To study this question, take a segment $<[1:0]; [0:1]> \subset \mathbb{C} \mathbb{P}^1_w$, connecting the poles, and lift its preimage under $\Psi$. The preimage  divides
 $\Sigma$ into $d$ connected  open parts $B_1, ..., B_d$ such that for every $B_i$  for each $p \in B_i$ one has $\Psi^{-1}(\Psi(p))\cap B_i \equiv p$. On the other hand it exists $B_i$ such that
the symplectic area $\int_{B_i} \omega \geq \frac{m}{d}$. 

Then the existence of a smooth loop $\tilde \gamma_{Ch}$ which satisfies the properties 1) and 2) above is equivalent to the following numerical condition
$k \frac{m}{d} > 1$. Indeed,  every smooth loop satsifies 1) must lie in a component $B_i$ therefore the symplectic area of the dics  $\delta$ is strictly less than $\frac{m}{d}$. But  to have the right period
$k \delta$ must be integer, thus $k \frac{m}{d} > 1$. 

This inequality has clear geometric meaning: both the numbers are topological intersection indexes, where  $k m = {\rm ind} (D_b \cap \Sigma)$ and $d = {\rm ind} (D^+ \cap \Sigma)$. Note
since every $D_i$ and $\Sigma$ are complex, the topological intersections are given just by the numbers of the intersection points. Let $\rho_i = {\rm ind} (D_i \cap \Sigma)$;
then the last inequality means $\sum_{i=1}^s \rho_i > \sum_{\lambda_i > 0} \lambda_i \rho_i$.  

Note that all the constructions depend on one combinatorial data: relation $\sum_i \lambda_i [D_i] = 0$, which corresponds to integer vector $(\lambda_1, ..., \lambda_s)$.
We can say that $(\lambda_1, ..., \lambda_s)$ corresponds to ``middle'' case if the inequality ${\rm ind} (D_b \cap \Sigma) > {\rm ind} (D^+ \cap \Sigma)$ holds.
For such ``middle'' pseudotoric structure we claim that it does exist ${\rm BS}_{can}$ non standard lagrangian torus of Chekanov type.

For the ``middle'' pseudotoric structure  one can find an approriate ${\rm BS}_{can}$ non standard torus $T(\gamma_{Ch}, 0)$ with periods
$(\frac{1}{k}, 0, ..., 0)$; possibly we have another non stadard tori with periods say $(\frac{l}{k}, 0, ..., 0), l> 1$, but   we can apply the arguments from the above about
the monotonicity condition only for the  first torus with $l=1$. For this case again we have a holomorphic section $\alpha_b \in H^0(M, K^{-1}_M)$ invariant under the Hamiltonian action of $X_{f_i}$,
and again we can choose a parametrization $\theta$ for convex smooth loop $\tilde \gamma_{Ch}$ such that  
$\frac{\partial}{\partial \theta} (p) = {\rm cos}(\theta) X_{f_n}(p) + {\rm sin} (\theta) {\rm I} (X_{f_n}(p))$, and establish that the corresponding degree equals 1 (not $l$).

Therefore for this ``middle'' case we can claim the existence of monotonic non standard lagrangian tori of Chekanov type.   

\section{Non toric case}

The constructions can be directly generalized to  non toric  but pseudotoric   case. Below we study two examples: 4 - dimensional quadric $Q$ and full flag $F^3$ in $\mathbb{C}^3$.
As it was established  both are pseudotoric, see [1].

{\bf The third example.} Take the hypersurface $Q = \{ z_0 z_1 + z_2 z_3 + z_4 z_5 = 0 \} \subset \mathbb{C} \mathbb{P}^5$ with homogenious coordinates $[z_0: ... : z_5]$.
The restriction of the standard Kahler form $\omega$ is taken as the integer symplectic form. Note that $(Q, \omega)$ is monotonic: the anticanonical bundle is ${\cal O}(4)|_Q$
therefore $k = 4$.

The pseudotoric structure is given by the data $(B, \Psi, f_1, f_2, f_3)$ where $\Psi: Q \backslash B \to \mathbb{C} \mathbb{P}^1_w$ is given by the formulas 
$$w_0 = z_0 z_1, \quad w_1 =
z_2 z_3, \quad w_2 = z_3 z_4, \quad w_0 + w_1 + w_2 = 0,
$$ 
therefore the base set consists of 8 projective planes defined by the condition $z_0 z_1 = z_2 z_3 = z_4 z_5 = 0$, and
real Morse functions $f_1, f_2, f_3$ has the following form
$$
f_i = \frac{a_i (\vert z_0 \vert^2 - \vert z_1 \vert^2) + b_i (\vert z_2 \vert^2 - \vert z_3 \vert^2) + c_i (\vert z_3 \vert^2 - \vert z_4 \vert^2)}{\sum_{i=0}^5 \vert z_i \vert^2},
$$
where $(a_i, b_i, c_i)$ is $(1, 0, -1)$ and the cyclic permutations of it. 
 
In contrast with the toric case the map $\Psi$ has three singular fibers: over points $p_1 = [1:-1: 0], p_2 = [1: 0: -1]$ and $p_3 = [0: 1: -1]$. As it was shown in [1], the choice of
a smooth loop $\gamma \subset \mathbb{C} \mathbb{P}^1_w \backslash \{p_i \}$ and regular values $c_i$ of $f_i$ defines a smooth lagrangian torus $T(\gamma, c_i) \subset Q$. We are interested
in the ``middle'' values of $f_i$,  below we will explain the reason why we consider the case $c_i = 0$. Smooth loops on the complement $\mathbb{C} \mathbb{P}^1_w \backslash \{p_i \}$ are distinguished 
by their topological types, and we say that  $T(\gamma, 0)$ is of the standard type if $\gamma$ is non contractible or that $T(\gamma, 0)$ is of Chekanov type if it is contractible.

Now our task is to calculate the periods for both types in dependence on smooth loops $\gamma$. To do this we must find a ``section'' of map $\Psi$ as it was done
in the toric case. However in the present case we must as well calculate ``toric'' periods for $T(\gamma, 0)$ which means the periods of
3 - torus $T^3_0 = \{ f_1 = f_2 = f_3 = 0 \} \cap \overline{\Psi^{-1}(p)}$ which belongs to a toric fiber over $p \in \mathbb{C} \mathbb{P}^1_w, p \neq p_i$.   
For every such $T^3$ we have a distinguished basis $(a_1, a_2, a_3)$ represented by loops of the form $[e^{it}: e^{-it}: z_2: z_3: z_4: z_5], [z_0: z_1: e^{it}: e^{-it}: z_4: z_5]$
and $[z_0: z_1: z_2: z_3: e^{it}: e^{-it}]$ correspondingly; clearly the periods are $(0, 0, 0) {\rm mod} \quad \mathbb{Z}$.

To find the 4th period for $T(\gamma, 0)$ take Riemann surface $\Sigma \subset Q$ given by the intersection $\Sigma = Q \cap V$ where $V$ is the projective plane $V = \{ z_0 = z_1, z_2 = z_3, z_4 = z_5 \}
\subset \mathbb{C} \mathbb{P}^5$. As in the toric case it was choosen such that $V \subset N_0$ where $N_0 = \{ f_1 = f_2 = f_3 =0 \}$. By the very construction $\Sigma$
is a plane conic, therefore its symplectic area is $m = 2$. It is ramified over $\mathbb{C} \mathbb{P}^1_w$ but the covering structure is more complicated than in the toric case: we have 4 sheets outside of
$p_i$, and over each $p_i$ the four branches are separated in two pairs so the preimage over each $p_i$ consists of 2 points (and there are exactly three possibility
for this separation, and all of them are realized over $p_1, p_2$ and $p_3$).

Our conic $\Sigma$ has 6 distinguished branching points $p_i^{\pm}$ such that $p_i^{\pm} = \Psi^{-1}(p_i)$; it has three distinguished symmetries given by rotations
$\tau_i$ around axes $<p_i^+, p_i^->$ on angle $\pi$. Every contractible smooth loop $\gamma_{Ch} \subset \mathbb{C} \mathbb{P}^1_w \backslash \{p_i \}$
is lifted to  4 copies of smooth loops $\tilde \gamma_{Ch}, \tau_i (\tilde \gamma_{Ch})$; on the other hand every smooth contractible loop $\tilde \gamma_{Ch}$ such that
$\tilde \gamma_{Ch} \cap \tau_i(\tilde \gamma_{Ch}) = \emptyset$ for each $i = 1, 2, 3$, defines a smooth loop $\gamma_{Ch}$ on the base therefore it defines a smooth lagrangian
torus $T(\gamma_{Ch}, 0)$ of Chekanov type. Of course, the intersection $T(\gamma_{Ch}, 0) \cup \Sigma$ consists of all $\tilde \gamma_{Ch}, \tau_i (\tilde \gamma_{Ch})$,
but the period is the same since the symmetry. 

Remove from the big circle joining points $p_i$ on $\mathbb{C} \mathbb{P}^1_w$ one segment (say, with ends at $p_1$ and $p_2$) and then lift this cutting up to $\Sigma$;
 it leads to the division of $\Sigma$ into four pieces $B_1, ..., B_4$  corresponding to octahedron with removed 4 ``equatorial''
edges. The symplectic area of each $B_i$ equals to $\frac{1}{2}$; at the same time a smooth loop $\tilde \gamma_{Ch}$ sitting inside of $B_i$ satisfies the property 1);
since the symplectic area of $B_i$ is $\frac{1}{2}$  it is possible to choose a smooth loop $\tilde \gamma_{Ch} \subset {\rm int} B_i$
which bounds a disc of symplectic area $\frac{1}{4}$ --- and it is enough for us to constract the desired non standard lagrangian torus satisfies ${\rm BS}_{can}$ condition.
Indeed, for the corresponding torus $T(\Psi(\tilde \gamma_{Ch}), 0) \subset Q$ the periods are $(\frac{1}{4}, 0, 0, 0)$ thus it is ${\rm BS}_{can}$.

 At the same time there are  standard type lagrangian tori, defined by loops $\tilde \gamma^i$ ``centered'' at the branching points $p^{\pm}_i \in \Sigma$. Consider the case
$p_1^+ \in \Sigma \subset Q$, all the resting cases are essentially the same. A smooth loop $\tilde \gamma^1$, centered at $p_1^+$ is suitable for us 
if it is projected by $\Psi$ to a smooth loop $\gamma^1 \subset \mathbb{C} \mathbb{P}^1_w$; and to get a ${\rm BS}_{can}$ lagrangian torus we need $\tilde \gamma^1$ bounds a 
disc of symplectic area $\frac{k}{4}$ with integer $k$. Since the places of points
$p_i^{\pm} \in \Sigma$ are at the vertices of octahedron so the distance between any two points is fixed, we easily deduce that such loop $\tilde \gamma^1$ exists;
moreover since the total symplectic are of $\Sigma$ equals 2 there are three possibilitites: $\tilde \gamma^1$ can cut  $\frac{1}{8}$, $\frac{1}{4}$ or $\frac{3}{4}$ of the total symplectic area,
and all the cases lead to smooth lagrangian tori which are distinguished by their periods. Denote as $\tilde \gamma^1_k$ the loop which cuts $\frac{k}{8}, k = 1, 2, 3,$ of the total symplectic area. 
Take the images $\gamma_k^1 = \Psi(\tilde \gamma^1_k)$ and construct the corresponding lagrangian tori
$T(\gamma_k^1, 0), k = 1, 2, 3$; all of them are ${\rm BS}_{can}$. It is clear that $T(\gamma_1^1, 0)$ and $T(\gamma_2^1)$ are not Hamiltonian isotopic since they have different periods.

We can repeat the construction choosing another $p_i$ as the ``center'' of the corresponding loop, therefore we have 9 different
${\rm BS}_{can}$ lagrangian tori of standard type;   thus a natural question arises --- which of them are  monotonic?  Of course, essentially
the question is about $T(\gamma_k^1, 0)$  since all the others possess the same properties. 

 Certain arguments hints the answer: our quadric $Q$ admits a toric degeneration, being included in the family $Q_t = \{z_0 z_1 + z_2 z_3 + t z_4 z_5 = 0 \}, t = [0; 1]$.
All these $Q_t$ are pesudotoric: the pseudtoric structure is given by the same functions $f_1, f_2, f_3$ and the same map $\Psi$ restricted on each $Q_t$ idividually;
as the base one has $\mathbb{C} \mathbb{P}^1_t = \{ w_0 + w_1 + t w_2 = 0 \} \subset \mathbb{C} \mathbb{P}^2$. At the limit $t=0$ for $Q_0$ we have that it is toric and singular:
it contains singular locus $L_{Sing} = \{z_0=z_1= z_2 = z_3 = 0\} \subset Q_0$ which is a projective line. Note that the base set for each $t$ is the same: $B = \{ z_0 z_1 =
z_2 z_3 = z_4 z_5 = 0 \}$, so $\Psi = \Psi_t: Q_t \backslash B \to \mathbb{C} \mathbb{P}^1_t$.

For the limit case $t = 0 $ the map $\Psi_0: Q_0 \backslash B \to \mathbb{C} \mathbb{P}^1_0$ admits just two singular fibers over the poles $[1: -1: 0]$ and $[0: 0: 1]$, being toric.
The section of $\Psi_0$ is again given by the intersection with the plane $V = \{z_0 = z_1, z_2 = z_3, z_4 = z_5 \}$; in this case the section $\Sigma_0$ is the union of two conjugated
projective lines, intersecting at the preimage $\Psi^{-1}([0: 0: 1])$. Geometrically under the limiting process $Q_t \mapsto Q_0$ the equator of $\Sigma_t$, containing $p_2^{\pm}$ and
$p_3^{\pm}$ is shrinked to point $[0: 0: 0: 0: 1: 1] \in \Sigma_0$. 

Under the deformation $t \mapsto 0$ the standard lagrangian tori $T(\gamma_k^1, 0) \subset Q_1 = Q$ are deformed to ``really'' standard lagrangian tori in a toric manifold
which can be described in the standard way. Take a smooth loop $\tilde \gamma_k$ centered at $p_1^+ = [1: 1: i: i: 0: 0] \in \{ z_0 + i z_2 = 0 \} \subset  \Sigma_0$
 which is invariant under the rotations around the axe $<p_1^+, p_1^->$ and such that it cuts disc of symplectic area $\frac{k}{4}$. Then the corresponding standard lagrangian tori
$T(\Psi_0(\tilde \gamma_k), 0) \subset Q_0$ are ${\rm BS}_{can}$.

But in the toric setup we can distinguish monotonic cases: since here we have the fair toric action generated by $f_1, f_2, f_3$ and the resting $f_4$ which is taken to 
preserve the section $\Sigma_0$, f.e. we can take
$$
f_4 = \frac{\vert z_4\vert^2 + \vert z_5 \vert^2}{\sum_{i=0}^5 \vert z_i \vert^2},
$$
it exists an invariant section $\alpha_b$ of the anti canonical bundle, whose zeros form the boundary divisor $D_b \subset Q_0$, and for each standard torus $T(\gamma_k, 0)$
we have an invarinat section $\delta_b$ of the determinant bundle, given by the wedge product of the Hamiltonian vector fields $X_{f_i}$. Therefore the natural pairing
$\alpha_b(\delta_b) \in \mathbb{C}^*$ is a constant function on $T(\gamma_k, 0) \subset Q_0$: its derivation along each $X_{f_i}$ must be trivial. 

It follows that the Maslov index for the loop $\tilde \gamma_k \subset \Sigma_0$ must be equal to the index of the intersection of the line $\{ z_0 + i z_2 = 0 \} \subset \Sigma_0$
and the part $D^+$ of the boundary divisor which is given in our case as $\overline{\Psi^{-1}([1: -1: 0])} \subset Q_0$ (note that $D^+$ does not contain the singular set
of $Q_0$). The index equals to 2, thus to be a loop on a monotonic lagrangian torus our loop must cut a disc of the symplectic area $\frac{1}{2}$ only. Consequently we have only one monotonic lagrangian torus
$T(\gamma_2^1, 0) \subset Q_0$.

Note that the universal Maslov class is stable with respect to continuous deformations; thus if we take the corresponding universal Maslov classes
for every standard lagrangian torus $T_t(\gamma_2, 0) \subset Q_t$ constructed using the same data on each $Q_t$ then if the universal Maslov class
for $T_0(\gamma_2, 0)$ is trivial it must be trivial for each $T_t(\gamma_2, 0), t \in [0; 1]$; and the opposite is also true.

These arguments lead to the following answer: for our  quadric $Q = Q_1$ the standard lagrangian tori $T(\gamma^i_2, 0) \subset Q, i = 1, 2, 3,$
are monotonic.

{\bf Remark.} Smooth standard lagrangian tori $T(\gamma^i_1, 0)$ and $T(\gamma_3^i, 0)$  in quadric $Q$ present for us quite interesting examples of ${\rm BS}_{can}$
lagrangian tori which carry non trivial universal Maslov classes.

{\bf The last example.} The full flag variety $F^3$ realized as ${\cal U} \subset \mathbb{C} \mathbb{P}^2 \times \mathbb{C} \mathbb{P}^2$, a divisor of type (1,1) given by the equation
$\sum_{i=0}^2 x_i y_i = 0$. Here $[x_0: x_1: x_2]$ and $[y_0: y_1: y_2]$ are homogenious coordinates on the first and the second direct summands correspondingly.
We consider the symplectic structure of type (1,1) given by the restriction to ${\cal U}$ of the product symplectic form $p_1^* \omega_x \oplus p_2^* \omega_y$.

The psedutoric structure constructed in [1] is given by two moment maps
$$
f_1 = \frac{\vert x_0 \vert^2 - \vert x_2 \vert^2}{\sum_{i=0}^2 \vert x_i \vert^2} - \frac{\vert y_0 \vert^2 - \vert y_2 \vert^2}{\sum_{i=0}^2 \vert y_i \vert^2}, \quad
f_2 =   \frac{\vert x_0 \vert^2 - \vert x_1 \vert^2}{\sum_{i=0}^2 \vert x_i \vert^2} - \frac{\vert y_0 \vert^2 - \vert y_1 \vert^2}{\sum_{i=0}^2 \vert y_i \vert^2};
$$
both the functions preserve by the Hamiltonian action the fibers of the map $\Psi: {\cal U} \backslash B \to \mathbb{C} \mathbb{P}^1_w \subset \mathbb{C} \mathbb{P}^2$
given by the equations $w_i = x_i y_i, \sum_{i=0}^2 w_i = 0$ where $[w_0: w_1: w_2]$ are the coordinates on the last $\mathbb{C} \mathbb{P}^2$. The base set $B$ consists
of six lines of the form $\{x_i = x_j = y_k =0 \}$ or $\{ x_i = y_j = y_k = 0 \}$ where $(i, j, k)$ is a permutation of $(0, 1, 2)$. 

The compactification of generic fiber $\overline{\Psi^{-1}(p)} = \Psi^{-1}(p) \cup B$ is isomorphic to del Pezzo surface $\mathbb{C} \mathbb{P}^2_3$ of degree 6; over three distinguished points
$p_1 = [1: -1: 0], p_2 = [1: 0: -1], p_3 = [0:1: -1] \in \mathbb{C} \mathbb{P}^1_w$ we have degenerated fibers isomorphic to 
two copies of del Pezzo surfaces  $\mathbb{C} \mathbb{P}^2_1 \cup \mathbb{C} \mathbb{P}^2_1$ both of degree 8; the details can be found in [1].

  The flag variety ${\cal U}$ is monotonic (being Fano variety), the monotonicity coefficient $k$ equals to 2: the anticanonical class $K^{-1}_{F^3}$ is presented by the restriction
${\cal O}(2,2)|_{{\cal U}}$ while our symplectic form represents $c_1({\cal O}(1,1)|_{{\cal U}})$.

Again, as it was done in the examples above, first we construct an approriate ``section'' of the map $\Psi$ on the common level set $f_1 =f_2 =0$
of our moment maps $f_1$ and $f_2$. In this example the suitable section is given by the following rational curve $\Sigma \subset {\cal U}$:
take the diagonal $\Delta = \mathbb{C} \mathbb{P}^2 = \{ x_i = y_i \} \subset \mathbb{C} \mathbb{P}^2 \times \mathbb{C} \mathbb{P}^2$ and intersect it
with our cycle, so $\Delta \cap {\cal U} = \Sigma$. It is clear that $f_i|_{\Sigma} = 0$; since in the diagonal $\Delta$ our curve $\Sigma$ is given
by the equation $\sum_{i=0}^2 x_i^2 = 0$ it is a plane conic, therefore it is rational. The symplectic area of $\Sigma$ equals to 2. 

Each regular fiber $\Psi^{-1}(p)$ intersects $\Sigma$ exactly at 4 points; we have ramification points at $p_1, p_2, p_3$ where 4 leaves are divided in pairs.
So geometrically the picture is essentially the same as in the previous example: the only difference appears when we calculate the ratio
between the symplectic area of $\Sigma$ (which is 2) and the monotonicity coefficient $k$ (which is 2, not 4 as in the quadric case). So the main difference changes the answer
to the question about ${\rm BS}_{can}$ non standard lagrangian tori. Indeed, as we have seen in the previous example, our $\Sigma$ can be divided into four
domains $B_1, ..., B_4$ such that a smooth loop $\tilde \gamma \subset {\rm Int} B_i$ defines the corresponding  non standard lagrangian torus. But in the previous case
the symplectic area of such a loop was restricted by $\leq \frac{1}{2}$ and since the monotonicity coefficient was 4 we could solve the existence problem; in contrast
now we can not do it --- a smooth loop $\tilde \gamma \subset B_i$ with the right period does not exist. 

Therefore the full flag variety $F^3$ does not admit ${\rm BS}_{can}$ non standard lagrangian tori with respect to the given pseudotoric structure.
It implies non existence of monotonic non standard lagrangian tori with respect to the given pseudotoric structure.

Perhaps if we consider another pseudotoric structure the answer will be different.

However our construction gives   monotonic lagrangian tori of standard type in $F^3$. Indeed, consider $\Sigma$ with 6 marked ramification points $p_i^{\pm}$.
Explicitly take $p_1^+ = [1: i: 0] \times [1: i: 0] \in \Sigma \subset {\cal U}$ and find a smooth loop $\tilde \gamma_{Cl}^1 \subset \Sigma$ centered at $p_1^+$ which is symmectric under the rotations
around the axe $<p_1^+, p_1^->$ and which cuts the disc of symplectic area $\frac{1}{2}$, --- such a loop exists since  the total area of $\Sigma$ is 2. For each point $p \in \tilde \gamma_{Cl}^1$
of this loop take the 2 - torus spanned by the toric action of the moment maps  $f_1$ and $f_2$. Then collecting these 2 - tori along $\tilde \gamma^1_{Cl}$ we get a smooth torus
$T(\gamma^1_{Cl}, 0) \subset F^3$. We claim that this torus is smooth and ${\rm BS}_{can}$. Indeed, since the projection $\gamma_{Cl}^1 = \Psi(\tilde \gamma^1_{Cl}) \subset
\mathbb{C} \mathbb{P}^1_w \backslash \{ p_i \}$ is a smooth non contractible loop, the resulting torus must be smooth; on the other hand it has periods $(\frac{1}{2}, 0, 0)$ thus it is
${\rm BS}_{can}$. 

Moreover, since the flag variety $F^3$ admits toric degeneration comaptible with our construction we  can repeat the arguments from the above and deduce that lagrangian
 torus $T(\gamma_{Cl}^i, 0)$ is monotonic.  

Indeed, as it was done for the quadric case above, consider the following deformation family ${\cal U}_t = \{ x_0 y_0 + x_1 y_1 + t x_2 y_2 = 0 \} \subset \mathbb{C} \mathbb{P}^2
\times \mathbb{C} \mathbb{P}^2, t \in [0; 1]$. For $t=1$ we have our given flag variety ${\cal U}_1 = {\cal U}$; note that the pseudotoric structure $(f_1, f_2, \Psi_t, B)$ is defined
over each element of our family, and at the limiting point $t = 0$ we get a pseudotoric structure on a singular toric variety ${\cal U}_0 = \{ x_0 y_0 + x_1 y_1 = 0 \}$: the singular set
consists of exactly one point $[0:0:1] \times [0:0:1]$, and the corresponding convex polytop $P_{{\cal U}_0}$ is the famous Gelfand - Zeytlin polytop with one 4 - valent vertex which corresponds to
the singular point.  

The map $\Psi_0: {\cal U}_0 \backslash B \to \{ w_0 + w_1 = 0 \} \subset \mathbb{C} \mathbb{P}^2$ has the same base set $B$ as for every ${\cal U}_t$ consists of 6 lines, but for $t=0$
there are 2 points $p_1, p_2$ underlying non regular fibers: $p_1$ is the same as for ${\cal U}$ with coordinates $[1: -1: 0]$ and the second point $p_2 = [0:0:1]$. Take the same diagonal
$\Delta$  and cut the section $\Sigma_0 \subset {\cal U}_0$ as $\Delta \cap {\cal U}_0$ (it is clear that condition $f_1|_{\Sigma_0} = f_2|_{\Sigma_0} = 0$ holds). But for $t=0$
our curve $\Sigma_0$ is reducible being the union of two projective lines with three marked points $p_1^{\pm} = [1: \pm i: 0] \times [1: \pm i: 0]$ and $p_2^+ = [0:0:1] \times [0:0:1]$ ---
the singular point of the ambient variety ${\cal U}_0$. 

Take the loop $\tilde \gamma_{Cl} \subset \Sigma_0$ centered at $p_1^+$ such that it is symmetric with respect to rotations around the axe $<p_1^+, p_2^+>$ and such that it bounds 
a disc of symplectic area $\frac{1}{2}$ (so this is an equatorial loop). Then the corresponding standard lagrangian torus $T_0(\Psi_0(\tilde \gamma_{Cl}), 0)$ is monotonic.
Again we can apply the deformation arguments and show that smooth lagrangian torus $T(\gamma_{Cl}^i, 0) \subset {\cal U}_1 = {\cal U}$ must be monotonic as well.

{\bf References:}

[1]  N. A. Tyurin, {\it  “Pseudotoric structures: Lagrangian submanifolds and Lagrangian fibrations”}, Russian Math. Surveys, 72:3 (2017), pp. 513- 546;

[2] N. A. Tyurin, {\it “Universal Maslov class of a Bohr–Sommerfeld Lagrangian embedding into a pseudo-Einstein manifold”}, Theoret. and Math. Phys., 150:2 (2007), 278–287;

[3] M. Audin, {\it ``Torus Action on Symplectic Manifolds''}, Progress in Mathematics, 93, Birkhauser, Basel, 2004.

\end{document}